\documentclass{amsproc}
\usepackage[T1]{fontenc}
\usepackage[latin1]{inputenc}
\input xy
\xyoption{all}
\usepackage{amsmath}
\usepackage{amsfonts}
\usepackage{amsthm}
\usepackage{amssymb}
\usepackage{amscd}

\newtheorem*{theorem*}{Theorem}

\numberwithin{equation}{section}
\renewcommand\a{\alpha}         

\renewcommand\d{\delta}
\newcommand\e{\varepsilon}

\newcommand\D{\Delta}
\newcommand\G{\Gamma}
\newcommand\f{\frac}

\newcommand{\Z}{{\mathbb{Z}}}

\newcommand{\C}{{\mathbb{C}}}

\newcommand{\Q}{{\mathbb{Q}}}

\renewcommand\Re{\text{Re~}}
\renewcommand\Im{\mbox{Im~}}
\renewcommand\i{^{-1}}
\renewcommand\({\left(}
\renewcommand\){\right)}

\newcommand{\ignore}[1]{}
\newcommand{\myignore}[1]{}
\newcommand{\mymyignore}[1]{}

\newcommand{\ttwo}[4]{
\(\begin{smallmatrix}{#1} & {#2}
\\ {#3} & {#4} \end{smallmatrix}\)}
\makeatletter
\def\imod#1{\allowbreak\mkern5mu({\operator@font mod}\,#1)}
\makeatother
\newcommand\smallf[2]{{\textstyle{\frac{#1}{#2}}}}
\newcommand\srel[2]{\begin{smallmatrix} {#1} \\ {#2} \end{smallmatrix}}

\newtheorem{thm}[equation]{Theorem}

\newtheorem{lem}[equation]{Lemma}

\newtheorem{prop}[equation]{Proposition}
\newtheorem{definition}[equation]{Definition}

\newcommand{\thmref}[1]{Theorem~\ref{#1}}

\newcommand{\lemref}[1]{Lemma~\ref{#1}}
\newcommand{\secref}[1]{Section~\ref{#1}}

\def\mc{\Mr}

\def\mc{{\mathcal M}}

\def\mod{\mathop{\rm mod}\nolimits}
\def\sgn{\mathop{\rm sgn}\nolimits}

\begin{document}

\title{A $p$-adic integral for the reciprocal of $L$-functions
}

\author{Stephen Gelbart}
\address{Nicki and J. Ira Harris Professorial Chair, Department of Mathematics,
Ziskind Building, Room 256, Weizmann Institute of Science, Rehovot 76100 Israel}
\email{steve.gelbart@weizmann.ac.il}

\author{Stephen D. Miller}
\address{Department of Mathematics,
Hill Center-Busch Campus,
Rutgers, The State University of New Jersey,
110 Frelinghuysen Rd,
Piscataway, NJ 08854-8019}
\email{miller@math.rutgers.edu}
\thanks{Miller was partially supported
by NSF grant DMS-0901594.}

\author{Alexei Panchishkin}
\address{Institut Fourier, UMR 5582 du CNRS,
Université de Grenoble I,
BP 74, 38402 Saint-Martin d'Hères, France }
\email{Alexei.Pantchichkine@ujf-grenoble.fr}

\author{Freydoon Shahidi}
\address{Department of Mathematics, Purdue University,
150 N. University Street, West Lafayette, IN 47907-2067}
\email{shahidi@math.purdue.edu}
\thanks{Shahidi was partially supported
by NSF grants DMS-0700280 and DMS-1162299.}

\date{December 13, 2012}

\dedicatory{Dedicated to the memory of Ilya Piatetski-Shapiro}


\begin{abstract}
We introduce an analog of part of the Langlands-Shahidi method to the $p$-adic setting,  constructing reciprocals of certain $p$-adic $L$-functions using the nonconstant terms of the Fourier expansions of Eisenstein series.  We carry out the method for the group $SL_2$, and give explicit $p$-adic measures whose Mellin transforms are reciprocals of Dirichlet $L$-functions.  The formulas for these measures involve Fourier coefficients of Eisenstein series, plus a delicately chosen multiple of Haar measure necessary for boundedness.

\end{abstract}

\maketitle

\section{Introduction and Motivation}\label{sec1}

                   The theory of
automorphic forms and $L$-functions is now understood in some general settings, especially the role played by the Langlands Program.
One approach towards this program has been through various methods of exploring
 $L$-functions of a single complex variable attached to automorphic forms.
These have included the method of integral representations (e.g., the methods of Riemann, Hecke, Tate, and Rankin-Selberg) as well as the Langlands-Shahidi method of studying the Fourier expansions of Eisenstein series.
%

One (long term) goal is to study \emph{$p$-adic} $L$-functions using these techniques. In particular, we are interested in a $p$-adic analysis of the Langlands-Shahidi method, which we now review in the complex setting \cite{La,Eulerproducts,Sha81,Sha88,Sha90}.  Motivated by his study of Selberg's theory of particular Eisenstein series, Langlands investigated Eisenstein series on  quite general reductive groups $G$, induced from automorphic representations $\pi$ of  smaller reductive groups $M$.
These series
 generalize
 \begin{equation}\label{Ezs}
    E(z,s) \ \  = \ \ \smallf 12 \sum_{\srel{(c,d)\neq (0,0)}{\gcd(c,d)=1}} \f{y^s}{|cz+d|^{2s}} \ \, , \ \ \ \Re{s}\,>\,1\,,
\end{equation}
in the setting $G=SL(2)$ and $M=GL(1)$.

In general, the  spectral theory of self-adjoint operators produces a meromorphic continuation of the Eisenstein series. In \cite{Eulerproducts}, Langlands showed that the constant terms in a Fourier expansion of these Eisenstein series inherited this meromorphic continuation and  were expressible as a ratio of  a product of automorphic $L$-functions. While this eventually led to the proof of the meromorphy of each individual L-function in the constant term \cite{Sha88,Sha90}, the proof of the functional equation and their holomorphy had to wait. In fact, many years later,  Shahidi computed  {\it non-constant} Fourier coefficients of those Eisenstein series  and proved the analytic continuation and functional equations of many $L$-functions.  In the setting of $G=SL(2),M=GL(1)$ from (\ref{Ezs}), the non-holomorphic Eisenstein series $E(z,s)$
has a Fourier expansion whose first Fourier coefficient is
\begin{equation}
\label{firstfouriercoeff}
2\,{y}^{1/2}\,K_{s-{1/2}}(2{\pi}y)\,\frac{1}{\Gamma(s)\,\zeta(2s)}\,;\end{equation}
see \cite{GeSha} for more details.
In the present discussion, the most important feature is the appearance of the Riemann $\zeta$-function in the denominator.

To summarize,  the analytic continuation and functional equation of the Riemann $\zeta$-function can be derived from the same properties of Eisenstein series.  In this paper we shall see how this $SL(2)$ phenomena carries over $p$-adically.

Of course, the Riemann $\zeta$-function has several different analytic constructions, most of which are based on Poisson summation (for example, Riemann's own construction as the Mellin transform of Jacobi's $\theta$-function).  Poisson summation is naturally a distributional statement (in the sense of Laurent Schwartz) about summing Dirac $\d$-functions at every integer,
\begin{equation}\label{psfdist}
    \sum_{n\,\in\,\Z}\d_{n}(x) \ \ = \ \ \sum_{n\,\in\,\Z}e^{2\pi i n x}\,.
\end{equation}
Thus the distribution $\sum_{n\,\neq\,0}e^{2\pi i nx}$ is equal to the constant $-1$ between integers, making the continuous periodic functions $\sum_{n\,\neq\,0}e^{2\pi i nx}n^{-k}$, $k\ge 2$, equal to polynomials on the open interval $0<x<1$. These polynomials are constant multiplies of the Bernoulli polynomials, and take the value $2\zeta(k)$ at $x=0$, for $k$ a positive even integer.
It was later understood that this underlying mathematics can be viewed as a construction of the Riemann $\zeta$-function as a $p$-adic Mellin transform of a $p$-adic measure, often called ``Mazur's measure''
\cite{Ku-Le,Iw-annals,MaSwD,Kob84} (more on this in \secref{sec:reviewofMazur} below.)

Serre \cite{Se73} gave another construction of the  $p$-adic $\zeta$-function using the constant terms of  Eisenstein series and $p$-adic modular forms, though his methods  overlap nontrivially (in particular, they also rely on Poisson summation). With this motivation, we present a $p$-adic measure based on the  {\em nonzero} Fourier modes of holomorphic Eisenstein series, whose Mellin transform is the reciprocal of the $p$-adic $\zeta$-function.  Hopefully this can later be extended to the larger class of $L$-functions that arise from the $(G,M)$-pairs that the Langlands-Shahidi method treats in general.  Our theorem reads as follows.

\begin{thm} \label{th1}  Assume that $p>2$ is a regular prime (i.e., $p$ does not divide the class number of $\Q(e^{2\pi i/p})$).
\par (a) There exists an explicitly defined $p$-adic measure $\mu^*$ on $\Z_p^*$ whose Mellin transform is the reciprocal  of the $p$-adic Riemann $\zeta$-function, i.e.,
\begin{equation}\label{Thm11a}
\int_{\Z_p^*} x^{k-1}\, d\mu^* \ \ = \ \
(1-p^{k-1})^{-1}\,\zeta (1-k)^{-1}
\end{equation}
for all even positive  integers $k$.
Moreover, the measure $\mu^*$ can be expressed
in terms of  the nonzero Fourier coefficients of  classical holomorphic
Eisenstein series, corrected by an additive nonzero multiple of $p$-adic Haar ``measure''.

\par (b) More generally,  for any nontrivial Dirichlet character\footnote{By definition, $\chi(p)=0$ even for the trivial character $\bmod p^\a$ when $\a>0$.  We also tacitly identify finite order characters on $\Z_p^*$ with {\em primitive} Dirichlet characters.} $\chi \mod {p^\a}$, $\a>0$, a similar equality holds, namely that
\begin{equation}\label{Thm11b}
\int_{\Z_p^*}\chi(x)\, x^{{k-1}}\,
d\mu^*  \ \ = \ \
L(1-{k}, \chi)^{-1}
\end{equation} for any  positive integer $k$ satisfying the parity condition $\chi(-1)=(-1)^k$.
\end{thm}

{\bf Remarks:}

1)  The regularity assumption on $p$ is crucial, for otherwise some of the $p$-adic $L$-functions in (\ref{Thm11a}-\ref{Thm11b}) definitely vanish, and such a measure $\mu^\star$ cannot exist.  Nevertheless, the Eisenstein construction and calculation we give applies to irregular primes, producing merely a $p$-adic distribution (instead of a measure).
It would be interesting, however, to identify the singular behavior for irregular $p$.  The continuation of the reciprocal of the $p$-adic Riemann $\zeta$-function recently came up completely independently in \cite{concon}, where the authors face some similar analytic issues.

 \

2)
Our method unfortunately says very little about the boundedness of the explicitly-defined distribution $\mu^*$, at least not directly.  We had initially hoped to find a direct proof  from properties of the Eisenstein series it is defined from, but
 instead rely on results about invertibility within the Iwasawa algebra (see section~\ref{sec:proofsec}). Thus the main contribution here is the explicit nature of the measure, in particular its derivation from Eisenstein series.
  On the other hand,
 should such a proof of the analytic continuation of reciprocal $L$-functions be discovered,  the same invertibility mechanism may then give a new proof of the boundedness of Mazur's measure in the regular case. There are deep analytic issues regarding the continuation of reciprocal $L$-functions (such as Leopoldt's conjecture on their nonvanishing at $s=1$), which demonstrates the subtlety of the issue.

\

3)   Mazur's measure is closely related to the trigonometric series
$\sum_{n\neq 0}e^{2\pi i n x} n^{-k}$ mentioned in the opening paragraph.  Likewise, our measure is closely related to $\sum_{n\neq 0}e^{2\pi i n x}\mu(|n|)n^{-k}$, where $\mu(\cdot)$ denotes the M\"obius $\mu$-function.  While the first sum is essentially a polynomial, the latter is analytically very difficult because of the behavior of its coefficients.  For example, $\sum_{n\neq 0}|n|^{-s}=2\zeta(s)$ has a meromorphic continuation to the complex plane with only a single, simple pole, whereas the  location of the poles of  $\sum_{n\neq 0}\mu(|n|)|n|^{-s}=2\zeta(s)\i$ depends on the Riemann hypothesis.

\

4)
Our proof uses neither  spectral theory nor the general theory of Eisenstein series, but these may perhaps be ingredients in generalizations to higher rank groups.  In fact, since we have looked at $SL(2)$ we have really used only classical predecessors to the Langlands-Shahidi method.
Our estimation is that obtaining an explicit formula for the reciprocal measure of  $p$-adic $L$-functions from Eisenstein series for higher rank groups (like the Langlands-Shahidi method does in the complex case) will be less difficult than showing   the boundedness necessary for the reciprocals to have $p$-adic analytic continuations.  As we mentioned in remark 2), this boundedness would give $p$-adic analytic continuations of those $L$-functions.  Since the $L$-functions treated by the Langlands-Shahidi method are more numerous than the  ones that can be $p$-adically analytically continued at present,  this approach might lead to new results in that area. That is the main motivation for this work, and the present paper should be viewed as an exercise that carries out the formal aspects in the simplest  case.   Since $p$-adic interpolation is not expected for all automorphic $L$-functions, additional algebraicity hypotheses will be necessary in higher rank (for example, that the archimedean component be of cohomological type).

\

5) Langlands wrote a 1987 letter to the first author about a conversation with Coates which speculated about which of the two mathematical methods of constructing complex $L$-functions would be most useful to $p$-adic $L$-functions.  Coates had remarked to Langlands that only the technique of integral representations had been useful $p$-adically. Unfortunately, this conversation, and this letter, were forgotten by everyone with the passing
of years. Only in 2010, when this paper was already begun, did the first author come across this letter while looking for something completely different. Langlands's letter, which goes on to describe the desirability of having concrete evidence that the Langlands-Shahidi method might be used $p$-adically, then served as extra encouragement for this project.

\subsection*{Acknowledgements}

We would like to extend our thanks to Ching-Li Chai, Haruzo Hida, Fabian Januszewski, Robert Langlands, Barry Mazur, Christopher Skinner, Jerrold Tunnell, and Eric Urban for their helpful conversations.   We wish to give particular thanks to Siegfried B\"ocherer for vetting some earlier ideas, and to John Coates for carefully explaining to us how the existence of $\mu^*$ follows from facts about the Iwasawa algebra.

We are  very grateful for the opportunity to participate in this volume in memory of Ilya Piatetski-Shapiro.  All four of us were collaborators of Ilya, and have deep admiration for his profound contributions and leadership in automorphic forms.  May his work and courage be remembered for generations.

\section{Mazur's  measure  and its real analytic interpretation}\label{sec:reviewofMazur}

One of the key motivations and ingredients in Theorem~\ref{th1} is  Mazur's measure
$\mu_{1,c}$ (defined below in (\ref{Mazurk})), which has the property that for any finite order character $\chi$ of $\Z_p$ and positive integer $k$ having the same parity as $\chi$
\begin{equation}\label{Mazurintegral}
    \int_{\Z_p^*}\chi(x)\,x^{k-1}\,d\mu_{1,c} \ \ = \ \ -\,(1-\chi(c)\i c^{-k})\,(1-\chi(p)p^{k-1})\,L(1-k,\chi)\,,
\end{equation}
where $c>1$ is a fixed integer relatively prime to $p$ which we specialize in section~\ref{sec:proofsec}.
In preparation for the Fourier series argument in section~\ref{sec2}, in this section we give a proof of this formula in similar terms.  Some of this material is motivated by a question of Mazur concerning how the  construction of the Riemann $\zeta$-function in  \cite{infiniteorder} relates to his $p$-adic measure.

We return to the setting of the opening paragraph, specifically   identity (\ref{psfdist}), and recall that the classical Bernoulli polynomials $B_k(x)$ are antiderivatives of either side
in the sense that
\begin{equation}\label{bernoulliantideriv}
    B_k(x) \ \  = \ \ -\,\f{k!}{(2\pi i)^k}\,\sum_{n\,\neq\,0}\f{e^{2\pi i n x}}{n^k} \ , \ \ \ \text{for} \ 0\,<\,x\,<\,1\, \ \text{and} \ \, k \,\ge \, 0\,.
\end{equation}
The above identity actually holds at the end points $x=0$ and $x=1$ if $k\ge 2$; by convention $B_1(0)=B_1=-1/2$.  Note that its values at $x\in \Q$ are the {\em additively twisted} $L$-functions which naturally arise as precursors to Dirichlet $L$-functions. For example, the value at $x=0$ expresses the $k$-th Bernoulli number as $B_k=B_k(0)=-(1+(-1)^k)k!(2\pi i)^{-k}\zeta(k)$.  The functional equation for the $\zeta$-function then gives the famous formula
\begin{equation}\label{Bernoulliformulaforzeta}
\zeta(1-k) \ \ = \ \ (-1)^{k-1}\,\f{B_k}{k}\,,
\end{equation}
for $k> 0$.

The fact that these special $L$-function values are related to polynomials  is crucial for this  theory, since it connects them to simpler algebraic expressions.
The $k$-th Bernoulli distribution is defined on $\Z_p$ in terms of the special values
\begin{multline}\label{muBk}
    \mu_{B,k}(a+p^m\Z_p) \ \  := \ \ p^{m(k-1)}\,B_k(\smallf{a}{p^m})\, \ \  \ \text{for}\ \ 0\,\le \,a\,<\,p^m  \ \ \text{and} \ \ k \, \ge \, 1 \,\\
\( \    =   \ \ -\,\f{k!}{(2\pi i)^k}\,p^{m(k-1)}\,\sum_{n\,\neq\,0}e^{2\pi i n a/p^m}n^{-k}  \ \ \ \text{if} \ \,  \ k \,\ge \, 2 \)\!.
\end{multline}
The presentation of this ($p$-adic) distribution in terms of the Fourier series (\ref{bernoulliantideriv}) demonstrates its additivity:

\begin{lem}\label{additivityprop}
    Let $a(n)$, $n\neq 0$, be a bi-infinite sequence such that $|a(n)|=O(n^{1-\d})$ for some $\d>0$, and define
\begin{equation}\label{additivityprop1}
    \nu(a+p^m\Z_p) \ \ = \ \ p^{m(k-1)}\,\sum_{n\,\neq\,0} a(n)\,e^{2\pi i n a/p^m}\,n^{-k}
\end{equation}
for a fixed prime number $p$ and integers $k\ge 2$, $m\ge 0$.
Then $\nu$ is additive in the sense that
\begin{equation}\label{additivityprop2}
    \nu(a+p^m\Z_p) \ \ = \ \ \sum_{b\,=\,0}^{p-1}\nu(a+bp^m+p^{m+1}\Z_p) \ \ \ \text{for all~~} a\,\in\,\Z\text{~and~}m\ge 0
\end{equation}
if and only if the sequence satisfies the property
\begin{equation}\label{additivityprop3}
    a(np) \ \ = \ \ a(n) \ \ \ \text{for all~~}n\,\in\,\Z_{\neq 0}\,.
\end{equation}
\end{lem}
\noindent
{\bf Remark:} the bound on the coefficients is used solely to ensure absolute convergence.
\begin{proof}
The righthand side of  (\ref{additivityprop2}) equals
\begin{equation}\label{additivityprop4}
p^{(m+1)(k-1)}\,\sum_{n\,\neq\,0} a(n)\,n^{-k}\,\sum_{b\,=\,0}^{p-1}e^{2\pi i n(\f{a}{p^{m+1}}+\f{b}{p})}\,.
\end{equation}
The sum $\sum_{b\imod p}e^{2\pi i nb/p}$ equals $p$ if $p|n$, and is zero otherwise.  Hence (\ref{additivityprop4}) equals
\begin{equation}\label{additivityprop5}
p^{(m+1)(k-1)}p\,p^{-k}\,\sum_{n\,\neq\,0} a(np)\,n^{-k}\,e^{2\pi i n a/p^{m}}\,.
\end{equation}
The difference between this and the lefthand side of  (\ref{additivityprop2})  is
\begin{equation}\label{additivityprop6}
p^{m(k-1)}\,\sum_{n\,\neq\,0} [a(np)\,-\,a(n)]\,n^{-k}\,e^{2\pi i n a/p^{m}}\,.
\end{equation}
Clearly (\ref{additivityprop3}) implies this vanishes.  Conversely, if an absolutely convergent (hence continuous) Fourier series vanishes on the dense set of  rational numbers having denominator   a power of $p$, its coefficients must all be zero.  Because of our assumptions on the size of the $a(n)$, this demonstrates that (\ref{additivityprop2}) implies (\ref{additivityprop3}).
\end{proof}
  It follows from Lemma~\ref{additivityprop} that the distributions $\mu_{B,k}$, which have $a(n)$ constant, are   additive.
Mazur's measures $\mu_{k,c}$ are defined as convolutions of $\mu_{B,k}$ with a $\d$-measure.  More precisely,
\begin{equation}\label{Mazurk}
    \mu_{k,c}(U) \ \ := \ \ \mu_{B,k}(U) \, - \,  c^{-k}\,\mu_{B,k}(c\,U)\ , \ \ \  \ U \subset\Z_p\ \, \text{compact open}
\end{equation}
 \cite[Chapter 2]{Kob84}.  Unlike the $\mu_{B,k}$, the $\mu_{k,c}$ are measures, that is, bounded distributions, and can be integrated against continuous functions such as polynomials.  Mazur proved the identity $\mu_{k,c}=k x^{k-1}\mu_{1,c}$.
 It follows from (\ref{Bernoulliformulaforzeta}) and  (\ref{muBk})  that the integral of $\mu_{k,c}$ over $\Z_p^*$ is
\begin{equation}\label{mazurscalculation}
\aligned
  \mu_{k,c}(\Z_p^*) \ \ & = \ \  (1-c^{-k})\,\mu_{B,k}(\Z_p^*) \ \  = \ \ (1-c^{-k})\,[\mu_{B,k}(\Z_p)-\mu_{B,k}(p\Z_p)] \\
& = \ \  (1-c^{-k})\,[B_k(0) - p^{k-1}B_k(0)] \ \ = \ \ (1-c^{-k})\,(1-p^{k-1})\,B_k(0) \\
& = \ \ (-1)^{k-1}\,k\,(1-c^{-k})\,(1-p^{k-1})\,\zeta(1-k)
\endaligned
\end{equation}
for $k\ge 1$.  This is assertion (\ref{Mazurintegral}) for the trivial character $\chi$.  For the nontrivial character computations we rely on the following:

\begin{lem}\label{generalMellincalc}
Assume the hypotheses of Lemma~\ref{additivityprop}, in particular the definition of $\nu$ in (\ref{additivityprop1}).  Then the integral of $\nu$ against a nontrivial Dirichlet character $\chi$ of conductor $q=p^\a$
is
\begin{equation}\label{generalMellincalc1}
    \int_{\Z_p^*}\chi(x)\,d\nu \ \ = \ \ q^{k-1}\,\tau_\chi\,\sum_{n\,\neq\,0}a(n)\,\overline{\chi(n)}\,n^{-k}\,,
\end{equation}
where
\begin{equation}\label{tauchidef}
    \tau_\chi \ \ := \ \ \sum_{r\imod {q}} \chi(r)\,e(\smallf{r}{q})\ , \, \ \ \tau_\chi \tau_{\overline{\chi}} \ = \ \chi(-1)\,q\,,
\end{equation}
is the Gauss sum for $\chi$.
\end{lem}
\begin{proof}
Indeed,
\begin{equation}\label{nuvschi}
    \aligned
    \int_{\Z_p^*}\chi(x)\,d\nu \ \ & = \ \ \sum_{\srel{a\,=\,1}{\gcd(a,p)\,=\,1}}^{q}
    \chi(a)\,\nu(a+q \Z_p) \\
    & = \ \ \sum_{\srel{a\,=\,1}{\gcd(a,p)\,=\,1}}^{q}  \chi(a)\,q^{k-1}\,
    \sum_{n\,\neq\,0} a(n)\,e^{2\pi i n a/q}\,n^{-k}\\
 & = \ \  q^{k-1}\,\sum_{n\,\neq\,0} a(n)\,n^{-k}\sum_{\srel{a\,=\,1}{\gcd(a,p)\,=\,1}}^{q}  \chi(a)\,e^{2\pi i n a/q} \\
 & = \ \ q^{k-1}\, \sum_{n\,\neq\,0} a(n)\,n^{-k}\,\overline{\chi(n)}\,\tau_\chi\,.
\endaligned
\end{equation}
Here we have used the fact that for a primitive character $\chi$ of modulus $q$,
\begin{equation}\label{gausssumfact}
    \sum_{\srel{a\,=\,1}{p\nmid a}}^q \chi(a)\,e^{2\pi i n a/q} \ \ = \ \ \left\{
                                                                            \begin{array}{ll}
                                                                              \chi(n)\i\,\tau_\chi, & p\nmid n\,, \\
                                                                              0, & p\mid n\,,
                                                                            \end{array}
                                                                          \right.
\end{equation}
i.e., equals $\overline{\chi(n)}\tau_\chi$ for all $n\in\Z$.
\end{proof}

Consequently, if $k\ge 2$ and $\chi$ is a nontrivial Dirichlet character of conductor $q$ and parity $\chi(-1)=(-1)^k$, we may use (\ref{Mazurk}) to compute
\begin{equation}\label{chivsmuka1}
\aligned
    \int_{\Z_p^*}\chi(x)\,d\mu_{k,c} \ \ & =  \ \
    (1-\chi(c)\i c^{-k}) \,\int_{\Z_p^*}\chi(x)\,d\mu_{B,k}\\
&  = \ \  -\,\f{k!}{(2\pi i)^k}\,(1-\chi(c)\i c^{-k}) \, q^{k-1}\,\tau_\chi\,\sum_{n\,\neq\,0}\overline{\chi(n)}\,n^{-k}
\\
    & = \ \ -\,(1-\chi(c)\i c^{-k})\,2\,\f{q^{k-1}k!}{(2\pi i)^k}\,\tau_\chi\,L(k,\chi\i)\,.
\endaligned
\end{equation}
Recall the functional equation of an $L$-function of a Dirichlet character of conductor $q=p^\a$ and parity $\chi(-1)=(-1)^\e$, $\e\in\{0,1\}$:
\begin{equation}\label{functionalequation}
    (\smallf{\pi}{q})^{s-1/2}\, \G(\smallf 12(1-s+\e)) \, L(1-s,\chi\i) \ \ = \ \ \f{i^\e q^{1/2}}{\tau_\chi} \G(\smallf12(s+\e))\,L(s,\chi)\,.
\end{equation}
It follows from this and standard $\G$-function identities that for $k\equiv \e\imod 2$
\begin{equation}\label{applyfe}
  2\,\f{q^{k-1}(k-1)!}{(2\pi i)^k}\,\tau_\chi\, L(k,\chi\i) \ \ = \ \  L(1-k,\chi)\,.
\end{equation}
Thus
\begin{equation}\label{chivsmuka2}
    \int_{\Z_p^*}\chi(x)\,d\mu_{k,c} \ \ = \ \ -\,k \, (1-c^{-k}\chi(c)\i)\,L(1-k,\chi)\,,
\end{equation}
which proves (\ref{Mazurintegral}) when $k\ge 2$.  The remaining situation of $k=1$ and  odd Dirichlet characters can be handled via a direct computation using the formula $B_1(x)=x-1/2$ for $0\le x<1$.

\vspace{.3cm}

{\bf Remark:}
Distributions whose Mellin transforms give other interesting Dirichlet series can of course be created using
Lemma~\ref{generalMellincalc}:~the more difficult aspect is showing that they are measures (i.e., bounded).  In our case of the reciprocal of the $\zeta$-function, one would take an appropriate multiple of $a(n)=\mu(n')$, where $\mu$ is the classical M\"obius function and $|n|$ factors uniquely as $n'$ times a power of $p$.  (We actually take a  different but related expression coming from Eisenstein series in the following section.)

Distributions whose $p$-adic Mellin transforms give $L$-functions (such as in \lemref{generalMellincalc} or \lemref{lem:mukstarintegrals}) can be formed from Fourier series which have arithmetically interesting coefficients $a(n)$.  Aside from the $q$-expansions of holomorphic Eisenstein series that we use in the following section,  another source of such Fourier series are automorphic distributions \cite{miller-schmid}.  For example, suppose that $\sum_{n\neq0}c(n)\sqrt{y}e^{2\pi i n x}K_0(2\pi |n|y)$ is a non-holomorphic Maass form corresponding to an even icosahedral Galois representation.  Its automorphic distribution is the Fourier series $\sum_{n\neq 0}c(n)e^{2\pi i nx}$ formed by removing the special function $\sqrt{y}K_0(2\pi |n|y)$.  Though this series converges only in the sense of distributions {\em ala} Laurent Schwartz, its antiderivatives $(2\pi i)^{-k} \sum_{n\neq 0}c(n)n^{-k}e^{2\pi i nx}$ are continuous functions for $k>0$.  The paper \cite{miller-schmid-jams} uses these distributions to obtain analytic continuations of complex $L$-functions.  It would be interesting if
this analogy to the Bernoulli polynomials that arose after integrating (\ref{psfdist}) proved useful for $p$-adic analytic continuations of higher degree automorphic $L$-functions.

%
%
%
%
%
%

\section{Fourier coefficients of classical Eisenstein series}\label{sec2}

In this section we review the Fourier expansions   of classical
Eisenstein series, and define $p$-adic distributions $\mu^*_k$ which will ultimately be used to define the
measure $\mu^*$ asserted in \thmref{th1}.  We begin with some motivating remarks about Eisenstein series for the congruence subgroups
\begin{equation}\label{G0NG1Ndef}
\aligned
    \G_0(N) \ \ &  = \ \ \{\,\ttwo abcd \,\in\,SL(2,\Z) \,|\, \ttwo abcd  \equiv \ttwo{\star}{\star}{0}{\star}\imod N \} \\
 \text{and~~~}  \G_1(N) \ \ &  = \ \ \{\,\ttwo abcd \,\in\,SL(2,\Z) \,|\, \ttwo abcd  \equiv \ttwo{1}{\star}{0}{1}\imod N \} \\
\endaligned
\end{equation}
of $SL(2,\Z)$,  $N>0$.  Given a Dirichlet character $\chi$ of modulus $N$ (not necessarily primitive) and an integer $k\ge 3$, the weight $k$ holomorphic Eisenstein series for $\G_0(N)$ transforming according to $\chi$ is defined by the formula
\begin{equation}\label{weightkholomchidef1}
    E_k(\chi,N) \ \ : = \ \ \sum_{\ttwo abcd \,\in\,\G_\infty\backslash \G_0(N)} \chi(d)\,(cz+d)^{-k}\,,
\end{equation}
where $\G_\infty=\{\ttwo 1n01|n\in\Z\}$.  Two matrices in $SL(2,\Z)$ are left-equivalent under $\G_\infty$ if and only if they have the same bottom row, and a tuple $(c,d)$ is the bottom row of a matrix in $\G_0(N)$ or $\G_1(N)$ precisely when $\gcd(c,d)=1$ and it satisfies the respective congruence condition in (\ref{G0NG1Ndef}).  Thus the sum in (\ref{weightkholomchidef1}) can be rewritten as
\begin{equation}\label{weightkholomchidef2}
\aligned
     E_k(\chi,N) \ \ &  = \ \ \sum_{\srel{(c,d)\,\in\,\Z^2}{\srel{\gcd(c,d)\,=\,1}{N|c}}} \chi(d)\,(cz+d)^{-k}
\\
& = \ \ \sum_{b\,\in\,(\Z/N\Z)^*} \chi(b) \sum_{\srel{(c,d)\,\in\,\Z^2}{\srel{\gcd(c,d)\,=\,1}{N|c,\,d\equiv b\imod N}}} (cz+d)^{-k} \\
& = \ \ \sum_{b\,\in\,(\Z/N\Z)^*} \chi(b) \sum_{\srel{(c,d)\,\in\,\Z^2}{\srel{\gcd(c,d)\,=\,1}{d\equiv b\imod N}}} (Ncz+d)^{-k}\,.
\endaligned
\end{equation}
In the last step we have changed variables $c\mapsto Nc$, and used the fact that for $d$ relatively prime to $N$ one has $\gcd(c,d)=\gcd(Nc,d)$.
When $N=1$ and $\chi$ is the trivial character,
this definition recovers the usual weight $k$ holomorphic Eisenstein series for $SL(2,\Z)$.

It is well-known that the nonzero Fourier coefficients of $E_k(\chi,N)$ involve the reciprocals of $L$-functions for the Dirichlet character $\chi$. The sum over $d$ in (\ref{weightkholomchidef2}) thus suggests a $p$-adic integral for this $L$-value involving a measure defined by the inner summation.  Accordingly, for any prime $p$ and integer $m\ge 0$ define a
map ${\mathcal E}_{k,p^m}$ from $\Z/p^m\Z$ to  ${\mathcal M}_k(\G_1(p^m))$ by the formula
\begin{equation}\label{curlyEdefnew}
    {\mathcal E}_{k,p^m}(b) \ \ := \ \ \sum_{\srel{(c,d)\,\in\,\Z^2}{\srel{\gcd(c,d)\,=\,1}{d\equiv b\imod {p^m}}}} (p^mcz+d)^{-k}\,.
\end{equation}
For example, ${\mathcal E}_{k,p^m}(1)$ is an Eisenstein series for $\G_1(p^m)$.

\begin{prop}
For $b$ coprime to $p$
the series $\mathcal E_{k,p^m}$ has the Fourier expansion
\begin{equation} \label{curlyEqexp}
\aligned
{\mathcal E}_{k,p^m}(b) \ \ = & \ \
 \sum_{\srel{r\,\in\,\Z}{r\,\equiv\,b\imod{p^m}}}{r^{-k}}
\sum_{\srel{n\,>\,0}{\srel{n\,\mid\,r}{p\,\nmid\,n}}}\mu(n)
 \\ & + \ \f{(-2\pi i)^k}{p^{mk}\,\G(k)}\,\sum_{\srel{t\,\in\, (\Z/p^m\Z)^*}{dd'\,>\,0}}  c_t(p^m)\,\sgn(d)\,d^{k-1}\,e^{2\pi idtb/p^m}\,e^{2\pi i dd' z}, \\
\endaligned
\end{equation}
where the coefficients $c_t(N)$ for $t$ relatively prime to $N$ are defined by the formula
\begin{equation}\label{ctdef}
c_t(N) \ \ := \ \ \displaystyle\sum_{\srel{n\,>\,0}{tn\equiv 1 \bmod {N}}} \frac {\mu (n)}{n^k}\,.
\end{equation}
\end{prop}
\begin{proof}

  For $k \ge 3$, $N\ge 1$, and $a,b \in {\mathbb Z}/N\Z$,
define
\begin{equation}\label{inv2.1}
{E}_{{k},N}(z;a,b)
\ \ = \ \ \sum_{\srel{(c,d)\,\in\,\Z^2-\{0,0\}}{(c,d)\equiv (a,b)\imod N}} (cz+d)^{-{k}}\,,
\end{equation}
which is in the space $\mc_k(\Gamma(N))$ of modular forms of weight $k$ for $\Gamma(N)$ (see \cite[p.271]{Miy}).  Using the identity
\begin{equation}\label{hurwitzidentity}
    \sum_{n \in \Z}(z+n)^{-{k}}  \ \  = \ \  \frac{(-2\pi i)^{k}}{\G(k)}\, \sum_{m \,>\,0}m^{{k}-1}e^{2\pi i m z} \ \, \ \ \ \ \Im(z) \,>\,0\,,
\end{equation}
 the Fourier expansion of $E_{k,N}$ can be directly computed as follows \cite[p.201]{He27}:
\begin{equation}\label{heckeexpansion}
\aligned
\quad {E}_{k,N} (z;a,b)  \ \ = & \ \
\delta (\smallf{a}{N}) \, \times  \sum_{\srel{d\,\in\,\Z_{\neq 0}}{d\equiv b\imod N}} d^{-k}
 \  \ + \ \\ &
 \  + \ \
{ (-2 \pi i)^{{k}}  \over
N^{{k}} \Gamma ({k})}       \   \times  \sum_{\scriptstyle dd'>0
 \atop \scriptstyle d' \equiv a \bmod N} {\rm sgn } (d)\,
d^{{k}-1}e^{2\pi i(\f{db}{N} +\f{dd' z}{N})}\,,
\endaligned
\end{equation}
in which
$\delta(x)$ is the characteristic function of the integers.

We shall also require a variant of (\ref{inv2.1}) that is summed over relatively prime pairs:~for $a,b\in\Z$ such that $\gcd(a,b,N)=1$, define
\begin{equation}\label{EkNstar}
{E}_{{k},N}^*(z;a,b) \ \ = \ \ \sum_{\srel{(c,d)\,\in\,\Z^2}{\srel{(c,d)\equiv (a,b)\imod N}{\gcd(c,d)\,=\,1}}} (cz+d)^{-{k}}\,.
\end{equation}
Using the M\"obius function $\mu(\cdot)$  and the assumption that $\gcd(b,p)=1$, it  can also  be expressed in terms of (\ref{inv2.1}) as
 \begin{equation}\label{moebius}
{ E}_{{k},N}^* (z;a,b)
\ \ = \ \ \sum_{t\,\in\,(\Z/N\Z)^*} c_t(N)\, { E}_{k,{N}} (z;ta,tb)
 \ \ \in  \ \  \mc_k(\Gamma(N))\, .
\end{equation}
We now insert (\ref{moebius}) into (\ref{curlyEdefnew}) to obtain
\begin{equation}\label{curlyEqexp2}
    {\mathcal E}_{k,p^m}(b) \ \ =  \ \
\sum_{a\,\in\,\Z/p^m\Z} E_{k,p^m}^*(p^mz;a,b)
\ \ = \ \ \sum_{\srel{a\,\in\, \Z/p^m\Z}{t\,\in\,
(\Z/p^m\Z)^*}} c_t(p^m)\,E_{k,p^m}(p^mz;ta,tb)\,.
\end{equation}
According to (\ref{heckeexpansion}) its constant term only involves the summand for $a\equiv 0\imod {p^m}$, and equals
\begin{equation}\label{curlyEqexp3}
\gathered
    \sum_{t\,\in\,(\Z/p^m\Z)^*} \sum_{\srel{n\,>\,0}{n\,\equiv\,t\i\imod{p^m}}} \smallf{\mu(n)}{n^k}\sum_{\srel{d\,\in\,\Z}{d\,\equiv\,tb\imod {p^m}}}\smallf{1}{d^{k}} \ \ = \\
= \ \
\sum_{\srel{n\,>\,0}{p\,\nmid\,n}}\smallf{\mu(n)}{n^k}\sum_{\srel{d\,\in\,\Z}{nd\,\equiv\,b\imod{p^m}}}
\smallf{1}{d^k}
\\
= \ \
 \sum_{\srel{r\,\in\,\Z}{r\,\equiv\,b\imod{p^m}}}{r^{-k}}
\sum_{\srel{n\,>\,0}{\srel{n\,\mid\,r}{p\,\nmid\,n}}}\mu(n)\,.
\endgathered
\end{equation}This is the first line on the right hand side of (\ref{curlyEqexp}).

For the nonconstant terms, we insert the second line of (\ref{heckeexpansion}) into (\ref{curlyEqexp2}), getting
\begin{equation}\label{curlyEqexp4}
  { (-2 \pi i)^{{k}}  \over
p^{{mk}}\, \Gamma ({k})}  \,  \sum_{\srel{a\,\in\, \Z/p^m\Z}{t\,\in\,
(\Z/p^m\Z)^*}} c_t(p^m)\,
 \sum_{\scriptstyle dd'>0
 \atop \scriptstyle d' \equiv ta \imod{p^m}} {\rm sgn } (d)\,
d^{{k}-1}e^{2\pi i(\f{dtb}{p^m} +dd' z)}\,,
\end{equation}
which agrees with the second line on the right hand side of (\ref{curlyEqexp}) after changing variables $a\mapsto at\i\imod{p^m}$.

\end{proof}

Let ${\mathcal C}_n$ denote the map which sends a holomorphic modular form to the coefficient of $e^{2\pi i n z}$ in its Fourier expansion.  The proposition states that
\begin{equation}\label{firstcoeff}
\aligned
    {\mathcal C}_{p^m}\!\({\mathcal E}_{k,p^m}(b)\) \ \
& = \ \
\f{(-2\pi i)^k}{p^{mk}\,\G(k)}\, \sum_{j\,=\,0}^m p^{j(k-1)}
\sum_{\srel{n\,\neq\,0}{p\,\nmid\,n}}\f{\mu(|n|)}{n^k}\,e^{2\pi i \bar{n}p^j b/p^m}\,,
\endaligned
\end{equation}
where $\bar{n}$ denotes the modular inverse of $n$ mod $p^m$.
We define a distribution $\widetilde{\mu}_k^*$ on $\Z_p^*$ by setting
\begin{equation}\label{mukstardef}
\aligned
    \widetilde{\mu}_k^*(b+p^m\Z_p) \ \ : & = \ \  \f14 \,{\mathcal C}_{p^m}\!\({\mathcal E}_{k,p^m}(b)\) \\
 & = \ \ \f{(-2\pi i)^k}{4 \,p^{mk}\,\G(k)}\,\sum_{j\,=\,0}^m p^{j(k-1)}
\sum_{\srel{n\,\neq\,0}{p\,\nmid\,n}}\f{\mu(|n|)}{n^k}\,e^{2\pi i \bar{n}p^j b/p^m} \,.
\endaligned
\end{equation}
That $\widetilde{\mu}_k^*$ is additive follows immediately from its definition as a sum over a congruence class in (\ref{curlyEdefnew}), or can also be checked directly {\em ala} Lemma~\ref{additivityprop}.

\begin{definition}\label{def:mukstar}
Define a distribution on $\Z_p^*$ by the formula
\begin{equation}\label{modified1}
   \mu_{k}^* \ \ := \ \ \widetilde{\mu}_k^* \ + \ c(k)\,\mu_{\text{Haar}}\,,
\end{equation}
where $\mu_{Haar}$ is Haar ``measure'' on $\Q_p$ and   $c(k) = \f{p^{2k-1}}{p^k-1}\f{(1-p^{k-1})\i}{\zeta(1-k)}$ if $k$ is even, and 0 otherwise.
In explicit terms for $k$ even
\begin{multline}\label{modified8}
      \mu_{k}^*(b+p^m\Z_p)  \ \ = \\ \f{(2\pi i)^k}{4\G(k)}\sum_{\srel{n\,\neq\,0}{p\nmid n}}\f{\mu(|n|)}{n^k}\(\sum_{j\,=\,0}^m  p^{j(k-1)-mk}e^{2\pi i \bar{n} p^j b/p^m}
-
\f{p^{-m}}{1-p^{1-k}}
\).
\end{multline}
\end{definition}

\noindent Without the ``correction factor'' $c(k)\mu_{Haar}$, $\mu_k^*$ would not be a bounded distribution.  Indeed, the boundedness of $\mu^*$ (shown in the next section) implies that   $\widetilde{\mu}_k^*(b+p^m\Z_p)$ should be equal to $-c(k)p^{-m}$ plus a bounded $p$-adic number (with a bound that depends on $k$).  At present we do not know if this ``divergence'' can be directly seen from  (\ref{mukstardef}).

\begin{lem}\label{lem:mukstarintegrals}
For any Dirichlet character $\chi$ whose conductor is a power of $p$ and integer $k\ge 3$ of the same parity,
\begin{equation}\label{modified3}
    \int_{\Z_p^*}\chi(x)\i\,d\mu_k^*(x) \ \ = \ \ \f{(1-\chi(p)p^{k-1})\i}{L(1-k,\chi)}\,.
\end{equation}
\end{lem}

\noindent {\bf Remark:} This is not only the key computational step in our argument, but also serves the important purpose of establishing that the measure values $\mu_k^*(b+p^m\Z_p)$ are algebraic numbers -- a fact which is not obvious from their definition.  Indeed, the right hand side of (\ref{modified3}) is of course algebraic, and thus so is the integral of $d\mu_k^*$ against any linear combination of finite order characters with algebraic coefficients -- in particular   the characteristic function of $b+p^m\Z_p$.

\begin{proof}
We begin with the case of the trivial character and $k\ge 4$ even.  Using the additivity of $\widetilde{\mu}_k^*$ the formula
\begin{equation}\label{trivintmukstar1}
\aligned
     \widetilde{\mu}_k^*(\Z_p^*) \ \ & = \ \ \sum_{b\,=\,1}^{p-1} \widetilde{\mu}_k^*(b+p\Z_p) \\
& = \ \  \f{(2\pi i)^k}{4 \,\G(k)} \sum_{\srel{n\,\neq\,0}{p\,\nmid\, n}} \f{\mu(|n|)}{n^k}\sum_{j\,=\,0}^1
p^{j(k-1)-k}\sum_{b\,=\,1}^{p-1}e(np^jb/p) \\
& = \ \ \f{(2\pi i)^k}{4 \,\G(k)} \sum_{\srel{n\,\neq\,0}{p\,\nmid\, n}} \f{\mu(|n|)}{n^k}\(-p^{-k}\,+\,
\f{p-1}{p}\) \\
& = \ \ \f{(2\pi i)^k}{2 \,\G(k)} \,\f{(1-p^{-k})\i}{\zeta(k)}\,
\f{p^k-p^{k-1}-1}{p^k}
\endaligned
\end{equation}
readily follows from (\ref{mukstardef}).
The functional equation of the Riemann $\zeta$ function,
\begin{equation}\label{zetafe}
    \f{1}{\zeta(k)} \ \ = \ \ \f{2\, \G(k)}{(2\pi i)^k \,\zeta(1-k)}\,,
\end{equation}
then implies
\begin{equation}\label{trivintmukstar2}
\aligned
     \widetilde{\mu}_k^*(\Z_p^*) \ \  & = \ \ \f{(1-p^{-k})\i}{\zeta(1-k)}\,
\f{p^k-p^{k-1}-1}{p^k} \\
& = \ \  \f{(1-p^{k-1})\i}{\zeta(1-k)}\,
\f{(1-p^{k-1})\,(p^k-p^{k-1}-1)}{p^k-1}
\,.
\endaligned
\end{equation}
The integral of the trivial character over $\Z_p^*$ is this plus $c(k)$ times $\f{p-1}{p}$, the Haar measure of $\Z_p^*$, which totals to give $\f{(1-p^{k-1})\i}{\zeta(1-k)}$ as claimed.

Next suppose that $\chi$ is a nontrivial character of conductor $q=p^m>1$ and that  $k\ge 3$ has the same parity.  Since $\chi$ is orthogonal to Haar measure,  its integral against $\mu_k^*$ is equal to
\begin{equation}\label{chiint}
\aligned
    \int_{\Z_p^*} \chi(x)\i \,d\widetilde{\mu}_k^*(x) \ \ & = \ \ \f{(-2\pi i)^k}{4 \,q^k\,\G(k)}\, \sum_{\srel{b\,=\,1}{p\nmid b}}^{q}\chi(b) \sum_{\srel{n\,\neq\,0}{p\,\nmid\,n}}\sum_{j\,=\,0}^m p^{j(k-1)}
\f{\mu(|n|)}{n^k}\,e^{2\pi i \overline{nb}p^j/p^m}
    \\
& = \ \  \f{(-2\pi i)^k}{4 \,q^k\,\G(k)}\,  \sum_{\srel{n\,\neq\,0}{p\,\nmid\,n}}
\f{\mu(|n|)}{n^k} \sum_{j\,=\,0}^m p^{j(k-1)} \sum_{\srel{b\,=\,1}{p\nmid b}}^{q}\chi(b) \,e^{2\pi i \overline{nb}p^j/p^m}
\\
& = \ \ \f{(-2\pi i)^k}{4 \,q^k\,\G(k)}\,  \sum_{\srel{n\,\neq\,0}{p\,\nmid\,n}}
\f{\mu(|n|)}{n^k} \,\chi(n)\i\,\tau_{\chi\i}
\\
& = \ \ 2\,\tau_{\chi\i}\,\f{(-2\pi i)^k}{4 \,q^k\,\G(k)}\, \f{1}{ L(k,\chi\i)}
\\
    & = \ \  L(1-k,\chi)\i
\endaligned
\end{equation}
(the terms with $j>0$ vanished above because of (\ref{gausssumfact})).
In the last step we used the formula
\begin{equation}\label{applyfeinverse}
    \f{(-2\pi i)^k}{q^k\,\G(k)}\f{\tau_{\chi\i}}{L(k,\chi\i)} \ \ = \ \    \f{2}{L(1-k,\chi)}\,,
\end{equation}
which is a restatement of  the functional equation (\ref{applyfe}) in light of  (\ref{tauchidef}).
\end{proof}

\section{Proof of \thmref{th1}}\label{sec:proofsec}

We begin with some preliminaries about measures on $\Z_p^*$ and its subgroup $\G=1+p\Z_p$, which we recall is isomorphic to  $\Z_p$ as a topological group for $p>2$.  Let $\omega(x)=\lim_{n\rightarrow\infty}x^{p^n}$ denote the Teichm\"uller character of $\Z_p^*$; its image is the set $\D$ of $(p-1)$-st roots of unity in $\Z_p^*$.  The
group $\Z_p^*$ decomposes as $\D\times \G$, with the
map $x\mapsto \omega(x)\i x$ furnishing the projection  onto $\G$.

 Let $\nu$ denote a measure on $\G$, extended to the rest of $\Z_p^*$ by the  relation
\begin{equation}\label{proofsec1}
\nu(aU) \  \ = \ \ \nu(U)
\end{equation}
for any $a\in \D$ and compact open $U\subset\G$.
Since each $\omega^i$ is a continuous function on $\Z_p^*$, the product $\omega^i\nu$ is a measure on $\Z_p^*$ satisfying the transformation law
\begin{equation}\label{proofsec2}
(\omega^i\nu)(aU) \ \ = \ \ \omega(a)^i\,(\omega^i\nu)(U)\ \ = \ \ a^i\,(\omega^i\nu)(U)\, , \ \ a \,\in\,\D\,.
\end{equation}
Every measure $\mu$ on $\Z_p^*$ can be decomposed as a sum of measures of the form $\omega^i\nu$:
\begin{equation}\label{proofsec3}
    \mu \ \ = \ \ \f{1}{p-1}\, \sum_{i\,=\,1}^{p-1}\omega^i\nu_i\,,
\end{equation}
where each $\nu_i$ satisfies (\ref{proofsec1}) (indeed, take $(\omega^i \nu_i)(U)=\sum_{a\in\D}\omega(a)^{-i}\mu(aU)$).

Suppose now that $\chi$ is a continuous homomorphism from $\Z_p^*$ to $\C_p^*$, and let $\chi_0$ denote its restriction to $\G$.  An arbitrary element $x\in \Z_p^*$ factors as $\omega(x)\cdot \omega(x)\i x \in \D\times \G$.  The values of $\chi(\omega(x))$ are $(p-1)$-st roots of unity in $\C_p^*$, and hence $\chi$'s restriction to $\D$ has the form $x\mapsto x^{-j}$ for some $j\imod{p-1}$ (i.e., $x^j\chi(x)=1$ for all  $x\in\D$).
Since $\chi$ transforms under $\D$ by $\omega^{-j}$, the integral of $\chi$ against $\mu$ of the form (\ref{proofsec3}) only involves the term for $i=j$, and equals
\begin{equation}\label{proofsec4}
    \int_{\Z_p^*}\chi(x)\,d\mu \ \ = \ \  \f{1}{p-1}\,\int_{\Z_p^*}    \chi(x)\,d(\omega^j\nu_j)        \ \  =  \ \  \int_\G \chi_0 \,d\nu_j\,,
\end{equation}
reducing the integration of the character $\chi$ on $\Z_p^*=\D\times \G$ to that of $\chi_0$ on $\G$.

We next make the assumption that $\mu$ is an odd measure on $\Z_p^*$, in the sense that $\mu(-U)=-\mu(U)$; equivalently, $\nu_i\equiv 0$ if $i$ is even.
 We furthermore suppose that
 the restriction of each $\nu_i$ to $\G\cong \Z_p$ for $i$ odd is a unit in the Iwasawa algebra, which in the context of  $\G$ is a convolution algebra with respect to multiplication.
In terms of the isomorphism with formal power series, the invertibility condition is that the constant term
\begin{equation}\label{invertcondition}
    (\omega^i\nu_i)(\G) \ \ = \ \ \sum_{a\,\in\,\D}\omega(a)^{-i}\,\mu(a\G) \ \ = \ \ \int_{\Z_p^*}\omega(x)^{-i}\,d\mu
\end{equation}
is a $p$-adic unit.
In terms of measures, this invertibility means for each odd value of  $i\imod {p-1}$ Iwasawa's isomorphism theorem \cite[p.97]{lang} guarantees the existence of an ``inverse'' measure $\nu_i\i$ such that $\nu_i\star \nu_i\i$ is equal to the $\d$-distribution $\d_1$ at the identity, or in other words
\begin{equation}\label{proofsec5a}
     \ \ \int_\G \int_\G f(xy)\,d\nu_i(x)\,d\nu_i\i(y) \ \ = \ \ f(1)
\end{equation}
for any continuous function $f:\G\rightarrow\C_p$.
In particular when $i$ is odd,
\begin{equation}\label{proofsec5b}
    \int_\G \chi_0\,d\nu_i\i \ \ =  \ \ \( \int_\G \chi_0\,d\nu_i\)\i.
\end{equation}
After extending each $\nu_i\i$ to measures $\omega^i\nu_i\i$ on $\Z_p^*$ as in
(\ref{proofsec1})-(\ref{proofsec2}), we define a measure on $\Z_p^*$ by the formula
\begin{equation}\label{proofsec6}
    \mu\i \ \ := \ \ \f{1}{(p-1)}\sum_{\srel{1\le i \le p-1}{i\text{~odd}}}\omega^i \nu_i\i\,.
\end{equation}
By the same reasoning as in (\ref{proofsec4}) and assuming further that $\chi$ is odd,
\begin{equation}\label{proofsec7}
    \int_{\Z_p}\chi(x)\,d\mu\i \ \  = \ \ \int_\G\chi_0\,d\nu_j\i
\ \ = \ \ \(  \int_\G \chi_0\,d\nu_j\)\i  = \ \ \(  \int_{\Z_p}\chi(x)\,d\mu \)\i,
\end{equation}
that is, the integrals of $\mu$ and $\mu\i$ against $\chi$ are reciprocals of each other.

Let us now specialize the above discussion to the particular case that $\mu$ equals Mazur's measure $\mu_{1,c}$, which we recall is an odd measure of $\Z_p^*$.    By (\ref{invertcondition}), the condition that $\nu_i$ be a unit in the Iwasawa algebra is equivalent to the integral $\int_{\Z_p^*}\omega^{-i}(x)d\mu_{1,c}$ being a $p$-adic unit.  This integral  is computed in the case $k=1$ of (\ref{Mazurintegral}) as
$-(1-\chi(c)\i c\i)L(0,\chi)$, with $\chi$ equal to the nontrivial, odd Dirichlet character $\omega^{-i}$.
At this point we assume, as we may, that $c$ is a primitive root $\mod p$ which is not congruent to $\omega(c)\imod {p^2}$.  These assumptions are made to ensure that prefactor $\omega(c)^i c\i-1$ is  $p$ times a $p$-adic unit if $i=1$, and a $p$-adic unit if $i>1$.

At the same time,
 $pL(0,\omega\i)=-\sum_{a=1}^{p-1}\omega\i(a)a\equiv -\sum_{a=1}^{p-1}1\equiv 1 \imod{p}$ \cite[p.32]{washington}.
 This verifies the invertibility of $\nu_1$ in the Iwasawa algebra. For $i>1$, the Kummer  congruences show that   $p$ divides $L(0,\omega^{-i})=L(0,\omega^{p-1-i})$ if and only if $p$ divides the Bernoulli number $B_{p-i}$ \cite[Cor.~5.15]{washington}.  Since $p$ is regular, this does not happen and these $\nu_i$ are also invertible in the Iwasawa algebra.

  Therefore we conclude from (\ref{proofsec7}) that the measure  $\mu_{1,c}\i$  defined by (\ref{proofsec6}) satisfies the property
\begin{equation}\label{mazurinverse1}
    \int_{\Z_p^*} x^{k-1} \,\chi(x)\,d\mu_{1,c}\i \ \ = \ \ -(1-\chi(c)\i c^{-k})\i \,(1-\chi(p)p^{k-1})\i \, L(1-k,\chi)\i
\end{equation}
for any Dirichlet character $\chi$ and nonnegative integer $k$ of the same parity.
For the same reason as in (\ref{Mazurk}), the regularization
\begin{equation}\label{mustardef}
    \mu^*(U) \ \  := \ \ - \,\mu_{1,c}\i(U)\, +\,\f{1}{c}\, \mu_{1,c}\i(cU)
\end{equation}
(which is obtained by convolution with  the $\d$-measure concentrated at $c$) satisfies
\begin{equation}\label{mazurinverse2}
    \int_{\Z_p^*} x^{k-1} \,\chi(x)\,d\mu^* \ \ = \ \ (1-\chi(p)p^{k-1})\i \, L(1-k,\chi)\i\,.
\end{equation}
This shows the existence part of Theorem~\ref{th1}.

To conclude we shall match $\mu^*$ to the distribution constructed from Eisenstein series in section~\ref{sec2}.
Since $\mu^*$ is a measure, so is $x^{k-1}\mu^*$ for any integer $k$.  In particular, it is a distribution.  Let $\bar{\mu}_k^*$ denote the distribution on $\Z_p^*$ defined by the formula $\bar{\mu}_k^*(U):=\mu_k^*(U\i)$, so that the integral in  (\ref{modified3}) is equal to $\int_{\Z_p^*}\chi(x)d\bar{\mu}_k^*$.
Comparing (\ref{modified3}) and  (\ref{mazurinverse2}) for $k
\ge 3$, we see that $\bar{\mu}_k^*-x^{k-1}\mu^*$ vanishes when integrated against any Dirichlet character.  Since these span the space of locally constant functions on $\Z_p^*$, the distribution $\bar{\mu}_k^*-x^{k-1}\mu^*$  must be identically zero, and hence $\bar{\mu}_k^*$ and $x^{k-1}\mu^*$ are equal as distributions.  Both are thus bounded, making $\bar{\mu}_k^*$ a measure.  (We again remark that it would be highly desirable to have a proof of the boundedness of $\bar{\mu}_k^*$ that only uses properties of Eisenstein series.)  Thus for each $k\ge 3$ we conclude that $x^{1-k}\bar{\mu}_k^*$ is a measure  coinciding with $\mu^*$, proving the theorem.

\bibliographystyle{plain}

\end{document}